\documentclass[10pt,a4paper]{article}
\usepackage[margin=1.1cm]{geometry}
\usepackage{amsmath,amsthm,amsfonts,amssymb,amscd,cite,graphicx}
\usepackage{tikz}
\usetikzlibrary{decorations.pathreplacing}
\usepackage{pgfplots} 
\usetikzlibrary{fit}

\usepackage{url}
\usepackage{titlesec}
\titleformat{\section}
{\normalfont\fontsize{12}{15}\bfseries}{\thesection}{1em.}{}

\allowdisplaybreaks[4]

\let\oldbibliography\thebibliography
\renewcommand{\thebibliography}[1]{%
  \oldbibliography{#1}%
  \setlength{\itemsep}{-2pt}%
}

\begin{document}

\baselineskip=0.20in

\makebox[\textwidth]{%
\hglue-15pt
\begin{minipage}{0.6cm}	
\vskip9pt
\end{minipage} \vspace{-\parskip}
\begin{minipage}[t]{6cm}
\end{minipage}
\hfill
\begin{minipage}[t]{6.5cm}
\end{minipage}}
\vskip36pt

\noindent
{\large \bf $\boldsymbol{S}$-Motzkin paths with catastrophes and air pockets (very early version)}\\

\noindent
Helmut Prodinger$^{1,}\footnote{Corresponding author (hproding@sun.ac.za)}$\\

\noindent
\footnotesize $^1${\it Helmut Prodinger,
	Department of Mathematical Sciences, Stellenbosch University,
	7602 Stellenbosch, South Africa, and
	NITheCS (National Institute for Theoretical and Computational Sciences),
	South Africa}\\

\noindent
 (\footnotesize Received: Day Month 202X. Received in revised form: Day Month 202X. Accepted: Day Month 202X. Published online: Day Month 202X.)\\

\setcounter{page}{1} \thispagestyle{empty}

\baselineskip=0.20in

\normalsize

 \begin{abstract}
 \noindent
So called $S$-Motzkin paths are combined the concepts `catastrophes' and `air pockets.
The enumeration is done by properly set up bivariate generating functions which can be extended using
the kernel method.
 \\[2mm]
 \end{abstract}

\baselineskip=0.20in

\section{Introduction}

 Dyck paths consist of up-steps $(1,1)$ and down-steps $(1,-1)$, start at the origin and do not go below the $x$-axis; they  appear in many texts, we just give one major reference,
\cite{stanley-book}. Typically, the path returns to the $x$-axis at the end, but we also consider the scenario of open paths, where the paths end at level $j$, say.
A popular variation of Dyck paths are Motzkin paths; the difference is just that now a horizontal step $(1,0)$ is also allowed.

In this paper, we concentrate on $S$-Motzkin paths, which is a subfamily of all Motzkin paths: all three steps (up, level, down) must appear $n$ times, and, ignoring the down-steps, the sequence is $(1,0)(1,1)(1,0)(1,1)(1,0)(1,1)\dots(1,0)(1,1)$. The follow figure shows how these paths are recognized: The two layers enforce that the flat and up steps are interlaced. Only paths that end in the origin are $S$-Motzkin but we consider all paths wherever they end.

Here is an example of such an $S$-Motzkin path with 15 steps:
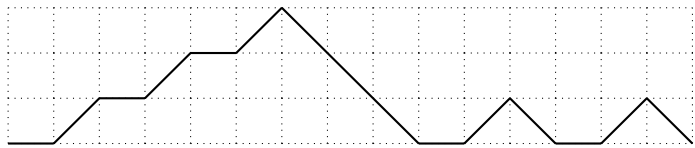
\begin{figure}[h]

	\begin{center}
		\begin{tikzpicture}[scale=0.6]

			\foreach \y in {0,1,...,3}
			\foreach \x in {1,2,...,15}
			{
				\draw[dotted] (\x-1,\y) -- (\x,\y);
			}
			\foreach \y in {0,1,...,2}
			\foreach \x in {0,1,2,...,15}
			{
				\draw[dotted] (\x,\y) -- (\x,\y+1);
			}
				\draw[thick] (0,0) to(1,0);	
							\draw[thick] (1,0) to(2,1);	
					\draw[thick] (2,1) to(3,1);	
	\draw[thick] (3,1) to(4,2);	
\draw[thick] (4,2) to(5,2);	
								\draw[thick] (5,2) to(6,3);	
											\draw[thick] (6,3) to(7,2);	
								\draw[thick] (7,2) to(8,1);	
							\draw[thick] (8,1) to(9,0);	
										\draw[thick] (9,0) to(10,0);	
															\draw[thick] (10,0) to(11,1);	
									\draw[thick] (11,1) to(12,0);	
											\draw[thick] (12,0) to(13,0);	
									\draw[thick] (13,0) to(14,1);	
											\draw[thick] (14,1) to(15,0);

		\end{tikzpicture}
	\end{center}
\end{figure}

Now we present a graph (automaton) to recognize exactly the $S$-Motzkin paths:
\begin{figure}[h]

	\begin{center}
		\begin{tikzpicture}[scale=1.5]

			\foreach \x in {0,1,2,3,4,5,6,7,8}
			{
				\draw (\x,0) circle (0.05cm);
				\fill (\x,0) circle (0.05cm);
				\draw (\x,-1) circle (0.05cm);
				\fill (\x,-1) circle (0.05cm);
			}

			\fill (0,0) circle (0.08cm);

			\foreach \x in {0,2,4,6}
			{
			}
			\foreach \x in {1,3,5,7}
			{
			}

			\foreach \x in {0,1,2,3,4,5,6,7}
			{
					\draw[thick, -latex] (\x+1,0) to  (\x,0);	
						\draw[thick, -latex] (\x+1,-1) to  (\x,-1);	
							\draw[thick, -latex] (\x,-1) to  (\x+1,0);	
				\node at  (\x+0.1,0.15){\tiny$\x$};
			}			
			
			\foreach \x in {0,1,2,3,4,5,6,7,8}
	{
							\draw[thick, -latex] (\x,0) to  (\x,-1);	
			}

			\node at  (8+0.1,0.15){\tiny$8$};

		\end{tikzpicture}
	\end{center}
	\caption{Graph to recognize $S$-Motzkin paths; they start and end at the special state (origin).}
\end{figure}
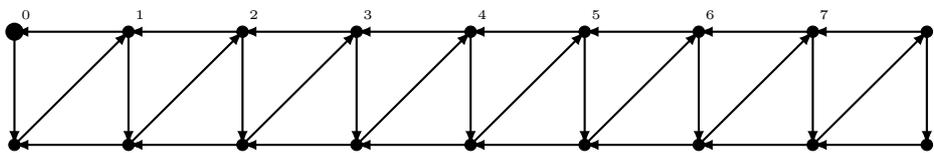

This subfamily of Motzkin paths originated from a question in a student competition; see \cite{PSW} and \cite{PUMA} for history and analysis.
In the following we will combine this family with \emph{catastrophes} and \emph{air pockets}, both originating in papers by Jean-Luc Baril and his team
\cite{BaryJIS}, \cite{baril-air}; the older paper by Banderier and Wallner \cite{BW} might be called the standard reference for lattice paths with catastrophes.
The very
recent papers \cite{Baryll, BaryJIS} contain some bijective aspects. Our own paper \cite{cata-Dyck}
investigates the situation in the context of skew Dyck paths.

Dyck (and other lattice) paths with catastrophes are  characterized by additional
steps (`catastrophes') that bring the path back to the $x$-axis in just one step from any level $j\ge2$.
For $S$-Motzkin paths the definition is similar, and the graphical description in    Figure \ref{purpel}
is easiest to understand;  the catastrophes are drawn in special colors.

In the last section, $S$-Motzkin paths with air pockets will be discussed. Briefly, down steps of any length are now allowed, but no
two down steps may follow each other.

\section{$S$-Motzkin paths with catastrophes}

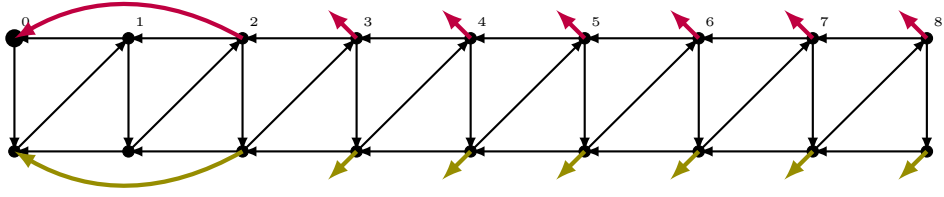
\begin{figure}[h]

	\begin{center}
		\begin{tikzpicture}[scale=1.5]

			\foreach \x in {0,1,2,3,4,5,6,7,8}
			{
				\draw (\x,0) circle (0.05cm);
				\fill (\x,0) circle (0.05cm);
				\draw (\x,-1) circle (0.05cm);
				\fill (\x,-1) circle (0.05cm);
			}

			\fill (0,0) circle (0.08cm);

			\foreach \x in {0,2,4,6}
			{
			}
			\foreach \x in {1,3,5,7}
			{
			}

			\foreach \x in {0,1,2,3,4,5,6,7}
			{
				\draw[thick, -latex] (\x+1,0) to  (\x,0);	
				\draw[thick, -latex] (\x+1,-1) to  (\x,-1);	
				\draw[thick, -latex] (\x,-1) to  (\x+1,0);	
				\node at  (\x+0.1,0.15){\tiny$\x$};
			}			
			
			\foreach \x in {0,1,2,3,4,5,6,7,8}
			{
				\draw[thick, -latex] (\x,0) to  (\x,-1);	
			}			
			
	\foreach \x in {3,4,5,6,7,8}
{
	\draw[ ultra thick, olive,-latex] (\x,0-1) -- (\x-0.25,-0.25-1);
	
}
\draw[ ultra thick, olive,-latex] (2,0-1) to[out=210,in=-30] (0, -1);

\draw[ ultra thick, purple,-latex] (2,0) to[out=150,in=30] (0, 0);
\foreach \x in {3,4,5,6,7,8}
{
	\draw[ ultra thick, purple,-latex] (\x,0) -- (\x-0.25,0.25);
	
}

			\node at  (8+0.1,0.15){\tiny$8$};

		\end{tikzpicture}
	\end{center}
	\caption{Graph to recognize $S$-Motzkin paths with catastrophes. Purple arrows lead to the initial state. Olive arrows lead to the level 0 state in the
	second layer.}
\label{purpel}
\end{figure}
In the sequel, we analyze the paths as in Fig.~\ref{purpel}.

We introduce generating functions $f_i=f_i(z)$, where the coefficient of $z^n$ counts the number of paths starting at the origin (=the big circle) 
and ending after $n$ steps at state $i$ (=level $i$) in the upper layer. 
The generating functions $g_i=g_i(z)$, where the coefficient of $z^n$ counts the number of paths starting at the origin (=the big circle) 
and ending after $n$ steps at state $i$ (=level $i$) in the lower layer. 

 The following recursions are easy to see:
\begin{align*}
	f_0&=1+z(f_1+f_2+f_3+f_4+\cdots),\\
	f_i&=zg_{i-1}+zf_{i+1}, \ i\ge1,\\
	g_0&=zf_0+z(g_1+g_2+g_3+g_4+\cdots),\\
	g_i&=zf_{i}+zg_{i+1}, \ i\ge1.
	\end{align*}
Since $f_0$ and $g_0$ are somewhat special, we leave them out for the moment and compute the other ones, $f_i$, $g_i$,  $i\ge1$.
Eventually we will solve the equations for $f_0$ and $g_0$, which will turn out to be just linear.
Therefore we introduce the bivariate generating functions
\begin{equation*}
	F(u)=F(u,z)=\sum_{i\ge1}u^{i-1}f_{i},\quad
	G(u)=G(u,z)=\sum_{i\ge1}u^{i-1}g_{i}
\end{equation*}
and we treat $f_0$ and $g_0$ as  parameters. Summing the recursions over all possible values of $i$,
\begin{equation*}
	F(u)=zg_0+zuG(u)+\frac zu[F(u)-f_1],\quad
	G(u)=zF(u)+\frac zu[G(u)-g_1].
\end{equation*}
Note that $f_1=F(0)$ and $g_1=G(0)$.
We compute
\begin{align*}
F(u)&=\frac{z(-u^2g_0+zug_0+uf_1-zf_1+zu^2g_1)}{z^2u^3-u^2+2zu-z^2},\\
G(u)&=\frac{z(-zu^2g_0+ug_1+zuf_1-zg_1)}{z^2u^3-u^2+2zu-z^2}.
\end{align*}
To factor the denominator, we set $u=zv$, and also $z^3=x=t(1-t)^2$ to get
\begin{equation*}
z^2(vt-1)(v^2t^2-2tv^2+vt+v^2-2v+1).
\end{equation*}
Therefore the three roots (expressed again in the variable $u$) are given by 
\begin{align*}
	u_1 =  \frac{z}{t}, \qquad u_2 = -z \frac{t - 2 + \sqrt{4t-3t^2}}{2(1-t)^2}, \qquad u_3 = -z \frac{t - 2 - \sqrt{4t - 3t^2}}{2(1-t)^2}
\end{align*}
and so
\begin{equation*}
z^2u^3-u^2+2zu-z^2=z^2(u-u_1)(u-u_2)(u-u_3).
\end{equation*}
These three roots appear already in \cite{PUMA}, were more details are provided.
Therefore
\begin{equation*}
F(u)=\frac{-u^2g_0+zug_0+uf_1-zf_1+zu^2g_1}{z(u-u_1)(u-u_2)(u-u_3)}\quad\text{and}\quad
G(u)=\frac{-zu^2g_0+ug_1+zuf_1-zg_1}{z(u-u_1)(u-u_2)(u-u_3)}.
\end{equation*}
Cancelling the bad factors $(u-u_2)(u-u_3)$ out, we get
\begin{equation*}
	F(u)=\frac{-g_0+zg_1}{z(u-u_1)}\quad\text{and}\quad
	G(u)=\frac{-g_0}{(u-u_1)}.
\end{equation*}
As a general remark, factors are bad if $\frac1{u-\overline u}$ has no power series expansion around $z=0$, $u=0$.
This is part of the kernel method, see \cite{Prodinger-kernel} for a user-friendly collection of examples. 
Plugging in $u=0$, we get
\begin{equation*}
	f_1=\frac{g_0-zg_1}{zu_1}\quad\text{and}\quad
	g_1=\frac{g_0}{u_1}=\frac{g_0t}{z}\quad\text{and thus}\quad f_1=g_0\frac{1-\frac{z}{u_1}}{zu_1}=
	g_0\frac{t(1-t)}{z^2}.
\end{equation*}
Now we can solve for $f_0$ and $g_0$:
\begin{align*}
	f_0&=1+z(f_1+f_2+f_3+f_4+\cdots)=1+zF(1)=1+\frac{-g_0+zg_1}{1-u_1}\\
	g_0&=zf_0+z(g_1+g_2+g_3+g_4+\cdots)=zf_0+\frac{-zg_0}{1-u_1},\\
\end{align*}
Therefore
\begin{equation*}
f_0=\frac{-t+z-zt}{-t+z-2zt+zt^2}\quad\text{and}\quad
g_0=\frac{z(z-t)}{-t+z-2zt+zt^2}.
\end{equation*}
Using the Lagrange inversion formula (or contour integration), we get the expansion
\begin{equation*}
t=\sum_{n\ge1}\frac1n\binom{3n-2}{n-1}x^n=\sum_{n\ge1}\frac1n\binom{3n-2}{n-1}z^{3n}.
\end{equation*}
\begin{align*}
f_0&=1+z^3+z^5+3 z^6+z^7+7 z^8+13 z^9+11 z^{10}+43 z^{11}+70 z^{12}+89 z^{13}+264 z^{14}+424 z^{15}+650 z^{16}+1657 z^{17}+\cdots\\
g_0&=z+2z^4+2z^6+7z^7+2z^8+15z^9+32z^{10}+23z^{11}+96z^{12}+174z^{13}+192z^{14}+604z^{15}+1048 z^{16}+1434 z^{17}+\cdots
\end{align*}
The coefficients are not `nice', in the sense that there are no simple expressions available for them. Consequently, $f_k$ and $g_k$ also do not have nice coefficients, although
the factor $\frac{1}{u-u_1}$ leads to nice coefficients, as can be seen from \cite{PUMA}.

Now we move to \textbf{asymptotics}.

As can be seen from the discussion in \cite{PSW}, the asymptotic enumeration of $S$-Motzkin paths is driven by a square-root type singulariy, as it often happens in the enumeration of trees and lattice paths:
\begin{equation*}
t\sim \frac13 - \frac2{3\sqrt3}\Bigl(1-\frac{27x}4\Bigr)^{1/2},
\end{equation*}
and the closest singularity (in $x$) is at $\frac{4}{27}$. Switching to the $z$-notation, as we have to in the context of catastrophes, we must look at the three 
roots closest to origin of modulus $\bigl(\frac{4}{27}\bigr)^{1/3}=0.5291336839$. Consequently, the exponential growth of $S$-Motzkin paths is given
by the reciprocal: $1.88988157485^n$. The exponent $n$ refers to the length $n$ of the $S$-Motzkin path. There are only paths when $n$ is divisible by 3, but
that is of no concern. 

For the case of catastrophes, we get a closer singularity. We need the dominant zero of the denominator $-t+z-2zt+zt^2$. A computer provides the value $\overline z=0.5248885986\dots$ and
the corresponding value $\overline{t}=0.2755080409\dots$\,. 
As we can see, the value is slightly smaller: $0.5248885986<0.5291336839$. Consequently this number leads to a \emph{simple} pole, and the exponential growth
is larger, as is not too surprising, considering the additional steps that are possible. The calculations are as follows:

We must expand $f_0$ and $g_0$ at the simple pole $z=\overline{z}$. First note that
\begin{equation*}
\frac t{dz}=\frac{dx}{dz}\frac {dt}{dx}\frac t{dt}=3z^2\frac1{(1-t)(1-3t)}\quad\text{and}\quad
\frac t{dz}\Big|_{z=\overline{z}, t=\overline{t}}=\frac{3\overline{z}^2}{(1-\overline{t})(1-3\overline{t})}=
\end{equation*}
\begin{align*}
-t+z-2zt+zt^2&\sim \frac d{dz}(-t+z-2zt+zt^2)\Big|_{z=\overline{z}}(z-\overline{z})
=\Big(-\frac t{dz}+1-2t-2z\frac t{dz}+t^2+2zt\frac t{dz}\Big)\Big|_{z=\overline{z}}(z-\overline{z})\\
&\sim-11.0530836206(z-\overline{z})\sim 21.0579609634(1-1.905166167z),
\end{align*}
and further
\begin{equation*}
f_0=\frac{-t+z-zt}{-t+z-2zt+zt^2}\sim 0.0049752931\frac{1}{1-1.905166167z}.
\end{equation*}
Since $f_0(z)$ is the generating function of $S$-Motzkin paths with catastrophes, we got an asymptotic equivalent of these numbers of length $n$ via
$[z^n]f_0(z)\sim0.0049752931(1.905166167)^n$. 

A similar computation leads to  $[z^n]g_0(z)\sim0.0062160344(1.905166167)^n$.
We continue with 
\begin{equation*}
	F(u)=\frac{-g_0+zg_1}{z(u-u_1)}=\frac{g_0(1-t)t}{z^2(1-u\frac tz)}
	\end{equation*}
and therefore
\begin{equation*}
[u^k]F(u)=\frac{(z-t)(1-t)}{-t+z-2zt+zt^2}\frac{t^{k+1}}{z^{k+1}}.
\end{equation*}
Similarly,
\begin{equation*}
	[u^k]G(u)=[u^k]\frac{t(z-t)}{-t+z-2zt+zt^2}\frac1{(1-\frac{ut}{z})}
	=\frac{t(z-t)}{-t+z-2zt+zt^2}\frac{t^{k}}{z^{k}}.
\end{equation*}
We note that $G(1)=\frac{zt}{-t+z-2zt+zt^2}$.
Now we move to partial $S$-Motzkin paths with arbitrary endpoint. In terms of generating functions, it just means $u:=1$, and we found the generating function
\begin{align*}
f_0(z)+F(1,z)+g_0(z)+G(1,z)&=\frac{ 1}{-t+z-2zt+zt^2}\Big[(-t+z-zt)+(1-t)t+z(z-t)+zt\Big]\\
&=\frac{ z+z^2-zt-t^2 }{-t+z-2zt+zt^2}\\
&=1+z+{z}^{2}+2{z}^{3}+3{z}^{4}+5{z}^{5}+10{z}^{6}+16{z}^{7}+
30{z}^{8}+58{z}^{9}+98{z}^{10}+189{z}^{11} +\cdots.
\end{align*}
 The asymptotic behaviour of the coefficients is also of the form $\textsf{const.}(1.905166167)^n$.


\section{Right-to-left S-Motzkin paths with catastrophes}

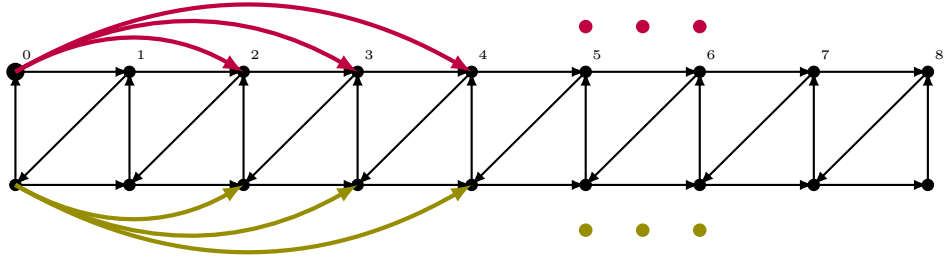
\begin{figure}[h]

	\begin{center}
		\begin{tikzpicture}[scale=1.5]

			\foreach \x in {0,1,2,3,4,5,6,7,8}
			{
				\draw (\x,0) circle (0.05cm);
				\fill (\x,0) circle (0.05cm);
				\draw (\x,-1) circle (0.05cm);
				\fill (\x,-1) circle (0.05cm);
			}

			\fill (0,0) circle (0.08cm);

			\foreach \x in {0,2,4,6}
			{
			}
			\foreach \x in {1,3,5,7}
			{
			}
			
\foreach \x in {0,1,2,3,4,5,6,7}
{
	\draw[thick, latex-] (\x+1,0) to  (\x,0);	
	\draw[thick, latex-] (\x+1,-1) to  (\x,-1);	
	\draw[thick, latex-] (\x,-1) to  (\x+1,0);	
	\node at  (\x+0.1,0.15){\tiny$\x$};
}

			\foreach \x in {0,1,2,3,4,5,6,7,8}
			{
				\draw[thick, latex-] (\x,0) to  (\x,-1);	
			}

			\draw[ ultra thick, olive,latex-] (2,0-1) to[out=210,in=-30] (0, -1);
			\draw[ ultra thick, olive,latex-] (3,0-1) to[out=210,in=-30] (0, -1);
			\draw[ ultra thick, olive,latex-] (4,0-1) to[out=210,in=-30] (0, -1);
			
			\draw[ ultra thick, purple,latex-] (2,0) to[out=150,in=30] (0, 0);
			\draw[ ultra thick, purple,latex-] (3,0) to[out=150,in=30] (0, 0);
			\draw[ ultra thick, purple,latex-] (4,0) to[out=150,in=30] (0, 0);

			\node at  (8+0.1,0.15){\tiny$8$};
			
				\foreach \x in {5,5.5,6}
			{
				\draw[ thick, purple] (\x,0+0.4) circle (0.05cm);
				\fill [ thick, purple](\x,0+0.4) circle (0.05cm);
				\draw[ thick, olive] (\x,-1.4) circle (0.05cm);
				\fill [ thick, olive](\x,-1.4) circle (0.05cm);
				
			}
			
		\end{tikzpicture}
	\end{center}
	\caption{Graph to recognize $S$-Motzkin paths with catastrophes from right-to-left.  }
	\label{purpel2}
\end{figure}
We use similar generating functions as before, namely $a_i(z)$ for the top layer, and $b_i(z)$ for the bottom layer. Then
\begin{align*}
a_0&=1+zb_0,\quad  a_1=zb_1+za_0, \quad a_i=zb_i+za_{i-1}+za_0,\ i\ge2,\\
b_0&=za_1,\quad b_1=zb_0+za_2,\quad b_i=za_{i+1}+zb_{i-1}+zb_0,\ i\ge2.
\end{align*}
As before, we introduce
\begin{equation*}
A(u)=\sum_{i\ge1}u^{i-1}a_i,\quad B(u)=\sum_{i\ge1}u^{i-1}b_i.
\end{equation*}
We compute
\begin{align*}
A(u)&=a_1+\sum_{i\ge2}u^{i-1}a_i=
a_1+z\sum_{i\ge2}u^{i-1}[b_i+a_{i-1}+a_0]\\
&=a_1+z\sum_{i\ge2}u^{i-1}b_i+z\sum_{i\ge2}u^{i-1}a_{i-1}+z\sum_{i\ge2}u^{i-1}a_0\\
&=a_1+zB(u)-zb_1+zuA(u)+\frac{zu}{1-u}a_0
\end{align*}
and
\begin{align*}
	B(u)&=b_1+\sum_{i\ge2}u^{i-1}b_i=b_1+z\sum_{i\ge2}u^{i-1}[a_{i+1}+b_{i-1}+b_0]\\
	& =b_1+z\sum_{i\ge2}u^{i-1}a_{i+1}+z\sum_{i\ge2}u^{i-1}b_{i-1}+z\sum_{i\ge2}u^{i-1}b_0\\
		& =b_1+\frac zu\sum_{i\ge0}u^{i}a_{i+1}-\frac zua_1-za_{2}+zuB(u)+\frac{zu}{1-u}b_0\\
			&  =\frac zuA(u)-\frac {z}ua_1 +zuB(u)+\frac{z}{1-u}b_0.
\end{align*}
We rewrite this system in the form
\begin{align*}
A(u)&=zuA(u)+zB(u)+\Phi(u),\quad \Phi(u)= a_1-zb_1+\frac{zu}{1-u}a_0 ,\\*
B(u)&=\frac zuA(u)+zuB(u)+\Psi(u),\quad \Psi(u)=-\frac {z}ua_1 +\frac{z}{1-u}b_0 .
\end{align*}
The solution is
\begin{equation*}
A(u)=\frac{z(zub_0-zu^2a_0-za_1+zua_1+ua_0)}{(z^2u^3-2zu^2+u-z^2)(1-u)},\quad
B(u)=\frac{z(-zu^2b_0-zu^2a_1+a_1u+b_0u+zua_1+za_0-a_1)}{(z^2u^3-2zu^2+u-z^2)(1-u)}.
\end{equation*}
Recall that if
\begin{align*}
	u_1 =  \frac{z}{t}, \qquad u_2 = -z \frac{t - 2 + \sqrt{4t-3t^2}}{2(1-t)^2}, \qquad u_3 = -z \frac{t - 2 - \sqrt{4t - 3t^2}}{2(1-t)^2};
\end{align*}
then
\begin{equation*}
z^2u^3-2zu^2+u-z^2=z^2\Big(u-\frac1{u_1}\Big)\Big(u-\frac1{u_2}\Big)\Big(u-\frac1{u_3}\Big);
\end{equation*}
this can be checked directly, compare also \cite{PUMA}. Since $\frac1{u_1}=\frac t{z}\sim z^2$, the factor $\big(u-\frac1{u_1}\big)$ is `bad' and must cancel out.
Plugging in $u=\frac tz$ into the numerators, we must get 0. 
Dividing out the factor $u-\tfrac tz$, the solutions now look  like
 \begin{equation*} 
 	A(u)=\frac{ z(zb_0+za_1-uza_0-a_0t+a_0)}{z^2\big(u-\frac1{u_2}\big)\big(u-\frac1{u_3}\big)(1-u)},\quad
 	B(u)=\frac{z(-zub_0+za_1-zua_1+a_1+b_0-tb_0-ta_1) }{z^2\big(u-\frac1{u_2}\big)\big(u-\frac1{u_3}\big)(1-u)}.
 \end{equation*}
Note that
\begin{equation*}
	\Big(u-\frac1{u_2}\Big)\Big(u-\frac1{u_3}\Big)=u^2+\frac{t-2}{z}u+\frac zt.
\end{equation*}
Now we can plug in $u=0$ to get
\begin{equation*}
	A(0)=a_1=\frac{ z(zb_0+za_1-a_0t+a_0)}{(1-t)^2},\quad\text{and}\quad
	a_1=\frac{z(zb_0+a_0-a_0t)}{-z^2+1-2t+t^2}.
\end{equation*}
Since $a_0=1+zb_0$, we get $a_1=\dfrac{z(2zb_0+1 -t-ztb_0)}{-z^2+1-2t+t^2}$.
We also get
\begin{equation*}
B(0)=b_1=\frac{z( za_1 +a_1+b_0-tb_0-ta_1) }{(1-t)^2}.
\end{equation*}
From $b_0=za_1$ we find
\begin{equation*}
b_0=\frac{z^2(zb_0+a_0-a_0t)}{-z^2+1-2t+t^2}=\frac{t(1-t)}{-t+z-2zt+zt^2} ,\quad\text{and}\quad
a_0=1+zb_0=f_0=\frac{-t+z-zt}{-t+z-2zt+zt^2}.
\end{equation*}
In principle, one could also write formul\ae\ for general $a_k$ and $b_k$, by using partial fraction decomposition. Since the results
look very complicated and do not provide extra insight, we refrain from giving such explicit results.

A brief comment about asymptotics: Since the denominators are the same as in the left-right instance, we get again an exponential behaviour, with the same
rate as before. The concept of open end does not make sense here since there are infinitely many such paths of a given length $n\ge1$.

While one might be tempted to attack the current question using some bijective tricks, it is worthwhile to note that our approach is very flexible,
and, e. g., subsets of the catastrophes may be considered, with little extra efforts.

\section{$S$-Motzkin paths and air pockets}

Now we move to another model popularized by Baril, namely introducing air pockets. These are maximal chains of downsteps, but this time only counted as one step. 
Since the main issue of $S$-Motzkin paths to keep the pattern flat, up, flat, up, flat, up, \dots alive, the downsteps live their own live, and we managed to 
construct a graph with 4 layers of states, describing all possible scenarios.
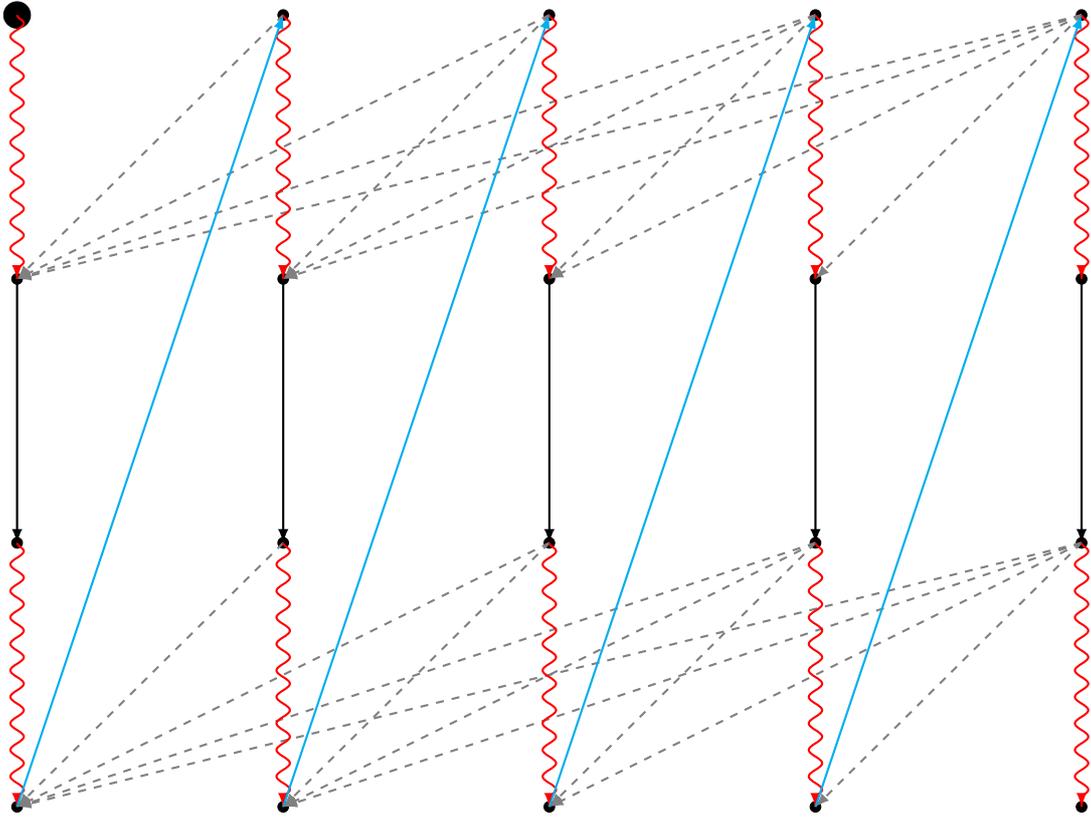
\begin{figure}[h]

	\begin{center}
		\begin{tikzpicture}[xscale=3.5,yscale=3.5]
			\draw (0,0) circle (0.05cm);
			\fill (0,0) circle (0.05cm);
			
			\foreach \x in {0,1,2,3,4}
			{
				\draw (\x,0) circle (0.02cm);
				\fill (\x,0) circle (0.02cm);
			}
			
			\foreach \x in {0,1,2,3,4 }
			{
				\draw (\x,-1) circle (0.02cm);
				\fill (\x,-1) circle (0.02cm);
			}
			
			\foreach \x in {0,1,2,3,4 }
			{
				\draw (\x,-2) circle (0.02cm);
				\fill (\x,-2) circle (0.02cm);
			}
			
			\foreach \x in {0,1,2,3,4 }
{
	\draw (\x,-3) circle (0.02cm);
	\fill (\x,-3) circle (0.02cm);
}

			\foreach \x in {0,1,2,3 }
			{

				
			}
			
			\draw[   thick, gray,dashed, -latex] (4,0) -- (3,-1);
			\draw[   thick, gray,dashed, -latex] (4,0) -- (2,-1);
			\draw[   thick, gray,dashed, -latex] (4,0) -- (1,-1);
			\draw[   thick, gray,dashed, -latex] (4,0) -- (0,-1);
			\draw[   thick, gray,dashed, -latex] (3,0) -- (2,-1);
			\draw[   thick, gray,dashed, -latex] (3,0) -- (1,-1);
			\draw[   thick, gray,dashed, -latex] (3,0) -- (0,-1);
			\draw[   thick, gray,dashed, -latex] (2,0) -- (1,-1);
			\draw[   thick, gray,dashed, -latex] (2,0) -- (0,-1);
			\draw[   thick,  gray,dashed, -latex] (1,0) -- (0,-1);

			\draw[   thick, gray,dashed, -latex] (4,0-2) -- (3,-1-2);
			\draw[   thick, gray,dashed, -latex] (4,0-2) -- (2,-1-2);
			\draw[   thick, gray,dashed, -latex] (4,0-2) -- (1,-1-2);
			\draw[   thick, gray,dashed, -latex] (4,0-2) -- (0,-1-2);
			\draw[   thick, gray,dashed, -latex] (3,0-2) -- (2,-1-2);
			\draw[   thick, gray,dashed, -latex] (3,0-2) -- (1,-1-2);
			\draw[   thick, gray,dashed, -latex] (3,0-2) -- (0,-1-2);
			\draw[   thick, gray,dashed, -latex] (2,0-2) -- (1,-1-2);
			\draw[   thick, gray,dashed, -latex] (2,0-2) -- (0,-1-2);
			\draw[   thick,  gray,dashed, -latex] (1,0-2) -- (0,-1-2);

			\foreach \x in {0,1,2,3,4 }
			{
				\draw [decorate,decoration=snake,thick	,red			 ] (\x,0) -- (\x,-1);
								\draw[   thick,red ,  -latex ] (\x,-1+0.05) -- (\x,-1);
																\draw[   thick ,  latex -] (\x,-2) -- (\x,-1);
																\draw [decorate,decoration=snake,thick	,red			 ] (\x,-2) -- (\x,-3);
																	\draw[   thick,red ,  -latex ] (\x,-3+0.05) -- (\x,-3);
			}
			
			\foreach \x in {0,1,2,3 }
			{
				\draw[   thick , cyan, latex -] (\x+1,0) -- (\x,-3);
			}

		\end{tikzpicture}
	\end{center}
	\caption{Four layers of states.}
	\label{threelayers2}
\end{figure}
Note that the wavy edges represent transitions without reading a symbol.
The generating functions for the four layers, reaching level $i$, can be read off from the diagram; note that the wavy edge is labelled by 1, not by $z$,
since there is no step done.
\begin{align*}
a_0&=1,\quad a_i=zd_{i-1},\ i\ge1,\quad
b_i=a_i+z\sum_{j>i}a_j,\\
c_i&=zb_i,\quad
d_i=c_i+z\sum_{j>i}c_j.
\end{align*}
The bivariate generating functions are $A(u)=\sum_{i\ge0}u^ia_i$, $B(u)=\sum_{i\ge0}u^ib_i$, etc. Summing the recursions over all values of $i$, we find the system
\begin{align*}
A(u)&=1+zuD(u),\quad
B(u)=A(u)+\frac{z}{1-u}[A(1)-A(u)],\\*
C(u)&=zB(u),\quad
D(u)=C(u)+\frac{z}{1-u}[C(1)-C(u)].
\end{align*}
The system can be reduced to two equations
\begin{align*}
	C(u)&=zA(u)+\frac{z^2}{1-u}[A(1)-A(u)],\\
	\frac{A(u)-1}{u}&=zC(u)+\frac{z^2}{1-u}[C(1)-C(u)];
\end{align*}
and so, by solving,
\begin{align*}
A(u)&={\frac {-{z}^{2}uC(1)+{z}^{3}{u}^{2}A(1)-1+2u+{z}^{4}uA(1)-{u}^{2}+{z}^{2}{u}^{2}C(1)-{z}^{3}uA(1)}{-1+2u-{u}^{2}+{
			z}^{2}u-2{z}^{2}{u}^{2}-2{z}^{3}u+{z}^{2}{u}^{3}+2{z}^{3}{u}^{2}
		+{z}^{4}u}},\\
	C(u)&={\frac { \left( {z}^{2}{u}^{2}C(1)-{u}^{2}-zu+2u+zA(1)u+{z
			}^{3}uC(1)-{z}^{2}uC(1)-1+z-zA(1) \right) z}{-1+2u-{u}^{
				2}+{z}^{2}u-2{z}^{2}{u}^{2}-2{z}^{3}u+{z}^{2}{u}^{3}+2{z}^{3}{u}
			^{2}+{z}^{4}u}}.
	\end{align*}
Plugging in $u=1$ gives the void equations $A(1)=A(1)$ and $C(1)=C(1)$. Therefore the denominator has to be investigated. We find that
$-1+2u-{u}^{2}+{	z}^{2}u-2{z}^{2}{u}^{2}-2{z}^{3}u+{z}^{2}{u}^{3}+2{z}^{3}{u}^{2}
+{z}^{4}u=z^2(u-\rho)(u-\sigma)(u-\tau)$; the explicit forms provided by Maple are useless, but fortunately gfun (in Maple) allows manipulations
with the relevant series:
\footnotesize
\begin{align*}
\rho&={z}^{-2}-2z-2{z}^{3}-{z}^{4}-2{z}^{5}-6{z}^{6}-4{z}^{7}-15
{z}^{8}-22{z}^{9}-33{z}^{10}-86{z}^{11}-115{z}^{12}-256{z}
^{13}-486{z}^{14}-804{z}^{15}-1783{z}^{16}-3074{z}^{17}\\&-6049
{z}^{18}-12104{z}^{19}-21902{z}^{20}-44918{z}^{21}-85235{z}^{
	22}-165124{z}^{23}-331137{z}^{24}-631740{z}^{25}-1261785{z}^{
	26}-2477694{z}^{27}+\cdots
\end{align*}
\normalsize
The other roots are ugly but we compute the simpler $(u-\sigma)(u-\tau)=u^2+Ku+L$ with
\begin{align*}
K&=-2{z}^{2}-2{z}^{5}-{z}^{6}-2{z}^{7}-6{z}^{8}-4{z}^{9}-15{
	z}^{10}-22{z}^{11}-33{z}^{12}-86{z}^{13}-115{z}^{14}-256{z}^
{15}-486{z}^{16}+\cdots,\\
L&={z}^{2}+2{z}^{5}+2{z}^{7}+5{z}^{8}+2{z}^{9}+14{z}^{10}+16
{z}^{11}+27{z}^{12}+74{z}^{13}+86{z}^{14}+222{z}^{15}+395{z}
^{16}+\cdots.
\end{align*}
So we must divide this term out from numerator and denominator.
Therefore
\begin{equation*}
A(u)=\frac{z^3A(1)+z^2C(1)-1}{z^2(u-\rho)},\quad\text{and}\quad
C(u)=\frac{z^2C(1)-1}{z(u-\rho)},
\end{equation*}
and by $u=1$ and solving,
\begin{equation*}
A(1)=\frac{-1+\rho}{z^2(-1+\rho+z)^2},\quad\text{and}\quad
C(1)=\frac{1}{z(-1+\rho+z)}.
\end{equation*}
Since these values are known, we find
\begin{equation*}
A(u)=\frac{\rho}{\rho-u}\quad\text{and}\quad C(u)=\frac{1-\rho}{z(1-\rho-z)(\rho-u)},
\end{equation*}
where the identity $1-2\rho+{\rho}^{2}-{z}^{2}\rho+2{z}^{2}{\rho}^{2}+2{z}^{3}\rho
-{z}^{2}{\rho}^{3}-2{z}^{3}{\rho}^{2}-{z}^{4}\rho=0$
was used for simplification. One sees $A(0)=1$, which is clear from combinatorial reasons.
Further
\begin{equation*}
	a_k=[u^k]A(u)=\rho^{-k}\quad\text{and}\quad c_k=[u^k]C(u)=\frac{1-\rho}{z(1-\rho-z)}\rho^{-k-1}.
\end{equation*}
The other quantities are then $b_k=\frac1zc_k$ and $d_k=\frac1za_{k+1}$, for any $k\ge0$.

We leave the analysis of this air pocket model from right to left, as well as other parameters, to the interested reader.
The factorization $(u-\rho^{-1})(u-\sigma^{-1})(u-\tau^{-1})$ will play a role here, and only one factor is bad, 
namely $(u-\rho^{-1})$.

It is possible to consider catastrophes and air pockets at the same time; we leave such considerations to enthousiastic 
young researchers.


\end{document}